\documentclass[a4wide]{article}
\pdfoutput=1

\usepackage[utf8]{inputenc}
\usepackage[T1]{fontenc}
\usepackage{graphicx,amssymb,amsmath}

\usepackage{xspace,bm}
\usepackage{array}

\usepackage{cite}

\newcommand{\abs}[1]{\lvert #1 \rvert}
\newcommand{\norm}[1]{\lVert #1 \rVert}

\usepackage{mathptmx}      % use Times fonts if available on your TeX system

\begin{document}
\pagestyle{myheadings}

\title{A Numerical Study of Newton Interpolation\\
       with Extremely High Degrees}

% all authors
\author{Michael Breuss, Friedemann Kemm and Oliver Vogel}

\date{}

% contact information of authors. Use one \contact command per person.
% Specify name(s), surname(s), postal address and e-mail.

% heading contains authors and title of the paper, please shorten the title if necessary
% (it will appear only at the heading of every even page)
\markboth{M.\ Breuss, F.\ Kemm and O.\ Vogel} {Newton Interpolation with Extremely High Degrees}

\maketitle

\begin{abstract}
In current textbooks the use of Chebyshev nodes with
Newton interpolation is advocated as the most efficient numerical
interpolation method in terms of approximation accuracy and computational
effort. However, we show numerically that the approximation quality
obtained by Newton interpolation with Fast Leja (FL) points is 
competitive to the use of Chebyshev nodes, even for extremely high 
degree interpolation. This is an experimental account of the 
analytic result that the limit distribution of FL points and 
Chebyshev nodes is the same when letting the number of points go to 
infinity. Since the FL construction is easy to perform and allows to add 
interpolation nodes on the fly in contrast to the use of Chebyshev 
nodes, our study suggests that Newton interpolation with FL points is 
currently the most efficient numerical technique for polynomial interpolation.   
Moreover, we give numerical evidence that any reasonable function can be 
approximated up to machine accuracy by Newton interpolation with FL points
if desired, which shows the potential of this method.
\end{abstract}

% here please insert two or more keywords appropriate to the topic of your paper

% here please insert one, two or more (5-digit) "Mathematics Subject Classification" codes appropriate to the topic of your paper
% you can find and choose them at the web page of the American Mathematical Society: http://www.ams.org/msc/msc2010.html
% \classification{00A01, 00M02, 00S03}
%\classification{65-05, 65D05, 97N50}

\section{Introduction}

Polynomial interpolation is one of the fundamental tasks in 
numerical analysis, see e.g.\ the textbooks 
\cite{atkinson89,gautschi97} for accounts on the subject. 
The goal of polynomial interpolation is to compute an approximation 
of an unknown function $f(x)$ 
over some one-dimensional real interval of interest $[a,b]$. 
This is to be done at hand of a set of 
given numbers $f\left( x_j \right)$
at \emph{interpolation nodes} $x_j$, $j=1, \ldots, n$, making use of 
a polynomial meeting exactly the data $(x_j, f\left( x_j \right))$.
A useful set of interpolation nodes is given by \emph{Chebyshev nodes}
which can be shown to minimize the oscillation of the node polynomial.
A popular format of the unique interpolation polynomial is the \emph{Newton form}
as its computation and evaluation can be done in an efficient way.
Still, the computation of high-degree interpolation polynomials for 
large numbers of given nodes is a non-trivial task since the process 
easily becomes numerically unstable, see e.g.\ \cite{tal-ezer91} for 
a useful discussion.

Let us turn to the approximation properties of an interpolation polynomial.
As Trefethen~\cite{trefethen-book} points out, in most cases the standard
polynomial interpolation is the most efficient numerical technique to
approximate continuous functions by polynomials. The difference
between approximation by Chebyshev polynomials and standard polynomial 
interpolation using Chebyshev nodes is rather small. In fact, it can be easily seen,
cf.~\cite[p.~389/399]{natanson}, that, if~\(\mathcal E_n\) denotes the
approximation error of the interpolation polynomial on Chebyshev nodes
with degree~\(n\), and~\(E_n\) denotes the error of the best
polynomial approximation, their ratio~\(\mathcal E_n / E_n\) is
bounded by~\(9+\frac{4}{\pi}\log n\). Hence, the loss in accuracy is
outweighed by the saving in terms of computational time. 
Moreover, one may expect that with increasing degree of an
interpolation polynomial in Chebyshev nodes, it will approach the underlying 
continuous function provided the computation can be done in a numerically stable 
and accurate way.

It is well-known that the location of the interpolation nodes $x_j$
influences the quality of the computed interpolant. However, not only the location 
but also the \emph{ordering} of the nodes can greatly affect the solution 
accuracy~\cite{calvettireichel03,reichel90}. 
While the Leja ordering of points is developed mainly for use in the complex domain (and in
general the ordering is not unique), also a version for use with real
intervals can be infered heuristically,
cf. \cite{baglama96,calvettireichel96,eisinberg06,higham90}. It puts
a fixed set of given points into a certain ordering, for instance it
can be applied as in this work at a set of predetermined Chebychev
nodes.

An undesirable property of the Leja ordering is that each time one
increases or decreases the number of nodes, the Leja ordering has to
be computed completely anew.  As a remedy to this, the \emph{Fast Leja (FL)}
points were proposed in \cite{baglama98}.  The procedure of
constructing the FL nodes involves selecting the location out of a
finite set of node candidates that is constructed in a simple way before 
the selection step. It also involves a fast way to determine the ordering during the
selection step. In the FL process, increasing the number of points
means to inherit the previously computed ones in the already
determined order. By these benefits, it seems of practical interest
to employ the FL points, especially with increasing degree of polynomial
interpolation. However, by the relatively simple construction
the question arises if the computations of the interpolation
polynomial can be done in a stable manner with no loss in approximation 
accuracy.

\paragraph{Our Contribution.}
In this paper we deal with the latter question by performing a thorough numerical
comparison of polynomial interpolation with FL points and the use of Chebyshev 
nodes sorted by the Leja ordering. In order to explain this proceeding it is important to 
note that the limit distribution of FL points -- i.e.\ the density of FL points as their number 
tends to infinity -- is the same as for Chebyshev nodes~\cite{reichel90}. 
Therefore, one can hope for the approximation properties of interpolation 
polynomials on FL points to be close to the approximation properties for 
interpolation on Chebyshev nodes, especially for high degrees of interpolation
polynomials. In fact, we give numerical evidence that any reasonable function can be
approximated up to machine accuracy by Newton interpolation on FL
points (at rather low numerical cost) at comparable quality as with 
Chebychev nodes. In doing this, we show that it is possible to give
a numerical analogon of the theoretical result that the limit
distribution of FL nodes tends to the distribution of Chebyshev nodes
in terms of approximation accuracy of corresponding polynomial interpolations.
By building exactly this bridge between theory and numerics we close the
corresponding slight gap in the literature.

\paragraph{Paper Organisation.}
In Section~\ref{sec:mathematical-basis} we give a brief account on the
aspects of the Newton polynomial interpolation problem, and we recall
the Leja ordering and the FL points.  This is followed by a
presentation of our computational study in
Section~\ref{sec:numer-exper}.  The paper is finished by a conclusion.

\section{Mathematical Basis}\label{sec:mathematical-basis}
%%%%%%%%%%%%%%%%%%%%%%%%%%%%%%

The purpose of this section is to briefly recall facts and methods as
employed in this work.

\subsection{Newton Polynomial Interpolation}\label{sec:newt-polyn-interp}
%%%%%%%%%%%%%%%%%%%%%%%%%%%%%%%%%%%%%%%%%%%%

Let us first recall the interpolation problem.
Given a discrete data point set $\{(x_j,f_j)\}$ for data points at 
position $x_j$ with function value $f_j\equiv f(x_j)$ and 
$j\in \{1, \ldots, n\}$, a polynomial $P_n(x)\in \Pi_n$ is sought. 
This polynomial should satisfy the interpolation condition~\(P_n(x_j) = f_j\)
% \begin{equation}
% \label{pol2}
% P_n(x_j) = f_j
% \end{equation}
in all data points. Gathering the 
interpolation conditions from all points, a linear system of
equations arises.

The \emph{Newton basis} 
consists of $n$ interpolation points
of the polynomials
\begin{equation}
\mathsf{n}_k(x) 
\; := \; 
\Pi_{i=1}^{k-1} \left( x - x_i \right) 
\qquad 
k=1,\ldots,n
\label{bvh-newton}
\end{equation}
The benefit of using this special basis is that the arising system
matrix is of lower triangular structure.  
Taking into account $n_0(x):=1$, it can easily be solved via
\begin{equation}
N(x_j)
\; := \;
\sum_{i=1}^{j} \mathsf{n}_{i}(x_j) a_{i} 
\; = \;
y_j 
\qquad 
\textrm{iterating as by}
\qquad
j = 1,\dots,n
\label{bvh-newtonpol}
\end{equation}
Moreover, for the evaluation of the Newton interpolation polynomial $N(x)$ 
the classical \emph{Horner scheme}~\cite{horner1819} can be used,
allowing for a stable, recursive $\mathcal{O}(n)$-algorithm 
to evaluate the polynomial.

\subsection{Chebyshev Nodes}\label{sec:chebyshev-nodes}
%%%%%%%%%%%%%%%%%%%%%%%%%%%%

So far, we did not explicitly mention how the interpolation points $x_j$ 
could be chosen. The most simple choice of equidistantly spaced points
is well known to lead to numerical instability, see for instance \cite{deboor01}. One 
way to cope with this is to consider the roots of 
Chebyshev polynomials 
\begin{align}
  T_m(t) &= \cos (m \arccos (t/b))\;, \qquad t\in [-b,b]\;,
\end{align}
with $m=0,1,\ldots$, where the roots are given by
\begin{eqnarray}\label{chebpoints}
  t_{k}^m := b \cos \left( \frac{(2k-1)\pi}{2m}\right) \qquad k = 1,
  2, \ldots, m
\end{eqnarray}
This choice of interpolation points is known to reduce the numerical
error problem~\cite{natanson,trefethen-book}. Let us note that a general, given interpolation
interval $[a,b]$ should be shifted and scaled as commented in Section~\ref{sec:numer-exper}.

For analytic functions, the error behaves like~\(\mathcal O(c^{-n})\) with some constant number $c$, 
which is rather fast. If~\(f\) is~\(\nu\)~times differentiable with~\(f^{(\nu)}\)
having bounded variation~\(V_\nu\), then the error is bounded
by~\(\frac{4V_\nu}{\pi\nu(n-\nu)^\nu}\). This implies a fast convergence
for~\(C^\infty\)-functions like the well known Runge function, which
cannot be approximated by interpolation on equidistant nodes. Even for
non-differentiable functions uniform convergence can be guaranteed if
the continuity module~\(\omega_f(\delta) = \sup_{\abs{x-y}\leq\delta}
\abs{f(x)-f(y)}\) decreases fast enough for~\(\delta \to 0\). For
functions with a finite number of discontinuities, we still have
point-wise convergence, but have to deal with the Gibbs phenomenon.
For a more detailed discussion we refer to~\cite{trefethen-book} and
the references therein.

\subsection{Leja Ordering of a Set of Points}\label{sec:leja-ordering}
%%%%%%%%%%%%%%%%%%%%%%%%%%

For $m$ given interpolation nodes in the set $S_m$, e.\,g.\ Chebyshev
nodes as defined above, one determines the \emph{Leja ordering} \cite{reichel90} of the
nodes $x_j$ via
\begin{eqnarray}
\abs{x_1}
%\left| x_1 \right| 
& := &
\max_{x \in S_m}
\abs{x}
%\left| x \right|
\label{bvh-leja1}\\
\underbrace{
\prod_{k=1}^{j-1}
\left|
x_{j}
-
x_k \right|}_{\textrm{'multiplicative distance'}}
& = &
\max_{l \textrm{ with }  j \leq l \leq m}
\prod_{k=1}^{j-1}
\left|
x_l - x_k
\right|
\, ,
\quad
2 \leq j \leq m
\, .
\label{bvh-leja2}
\end{eqnarray}
In order to implement this ordering, one simply starts with the 
largest point in the set of points. Then, one selects
the point that has the largest distance to this point. 
Usually, these two points will be the endpoints
of the interpolation interval $[a,b]$. Of the remaining 
nodes in the set to be sorted, one tests all for 
their multiplicative distances to the points selected so far.
Then one chooses the one with the largest multiplicative distance as the next
in the order. This process is repeated until all points 
are in order. 

Note, that sorting the points has a quadratic computational complexity. 
For Newton Interpolation, it pays off anyway since the Leja ordering considerably decreases rounding
errors compared to the naive ordering~\(x_0<x_1<\dots
<x_{n-1}<x_n\). While for the latter, Newton interpolation would fail
already for some~\(n<100\) even on Chebyshev nodes, for the Leja
ordering of Chebyshev nodes, as will be seen in
Section~\ref{sec:numer-exper}, there is no visible influence of
rounding errors at all.

\subsection{Fast Leja (FL) Points}\label{sec:fast-leja-fl}
%%%%%%%%%%%%%%%%%%%%%%%%%%%%%

FL points \cite{baglama98} represent an alternative to a Leja ordering
of Chebyshev nodes.  The algorithm to compute the FL points can be
described in a simple way.  The basic idea is to start with a small
set of candidate points~\(S^{(0)}\), which is filled with new
candidates in each iteration of the process.  For polynomial
Interpolation on the interval~\([a,b]\), one starts
with~\(S_2=S^{(0)}= \{a,b\}\). Then, one chooses the midpoint between
these two points as the first candidate point resulting in~\(S^{(1)} =
\{a,b,\frac{a+b}{2}\}\).  Now, one recursively chooses the candidate
point with the largest multiplicative distance of all points
in~\(S^{(n-1)}\setminus S_n\) as the next FL point~\(x_{n+1}\)
resulting in a new of nodes~\(S_{n+1}\). Again, one adds the midpoints
between this point and its neighbours in~\(S_{n+1}\) to the set of
candidate points, now~\(S^{(n)}\). Each~\(S_n\) is Leja ordered by
construction. A working Matlab code for this process and details on
the method can be found in \cite{baglama98}. 

A major advantage of FL points is the low computation cost and the
fact that they can be computed beforehand and stored for some large
enough~\(N\). Now every~\(S_n\) with~\(n\leq N\) just consists of the
first~\(n\) entries of~\(S_N\). This means also that, for Newton
interpolation, the addition of more nodes is rather simple, low cost,
and not bound to a certain number of additional points.

It is important to note that the limit
distribution, i.\,e.\ the density of FL points for~\(n\to\infty\) is
the same as for Chebyshev nodes~\cite{baglama98}. Therefore, one can hope for the
approximation properties of interpolation polynomials on FL points to
be close to the approximation properties for interpolation on
Chebyshev nodes.

%%%%%%%%%%%%%%%%%%%%%%%%%%%%%%%%%%%%%%%%%%%%%%%%%%%%%%%%%%%%%%%%%%%%%%

\section{Numerical Experiments}\label{sec:numer-exper}

In this section, we will demonstrate that very high degree polynomial
interpolation can be done using the numerical techniques presented earlier
in this paper.

\subsection{Technical Remarks}
For the evaluation of Newton interpolation polynomials of very high
degree the size and location of the computational domain may be point of concern.
It can be shown that a scaling of the $x$-domain to the length four
can be favorable. Moreover, the interpolation interval 
should ideally be identical to the interval $[-2,2]$, cf.\ \cite{tal-ezer91}.
Accordingly, we will employ such a shifting and scaling in our experiments.  

Also, it is necessary to employ
high-precision datatypes. Especially in the calculation of the Fast
Leja points and the Leja ordering of the Chebyshev nodes high
precision is needed in order to get the products in
equation~(\ref{bvh-leja2}) right. Therefore, we use for that task \emph{long
double}~(16~byte), and for the interpolation itself \emph{double}~(8~byte).

\subsection{The Interpolated Functions}
The numerical experiments are done on four different functions, all on the interval $[-2,2]$:
\begin{enumerate}
\item 
The Runge function \cite{runge1901}
\begin{equation}
f_1(x)=\frac 1 {1+6.25x^2}
\end{equation}
which is a classic example for a function which cannot be approximated
by interpolation on equidistant nodes. But, as already mentioned in
Section~\ref{sec:chebyshev-nodes}, it can be approximated at least by
interpolation on Chebyshev nodes. The purpose of this example is to
compare the convergence rate of interpolation on Chebyshev and Fast
Leja nodes. 
\item 
The Heaviside function
\begin{equation}
f_2(x)=\begin{cases} 1 &x>0\\0 & x\leq 0\, ,\end{cases}
\end{equation}
which has a discontinuity at the origin. Thus, we have to face the
Gibbs phenomenon and expect slow convergence near the
discontinuity. The same is true for the next example:
\item 
A sawtooth function
\begin{equation}
f_3(x)=x-[x]
\end{equation}
with several discontinuities.
\item
\begin{equation}
f_4(x)=\sqrt{|x|}
\end{equation}
which is not differentiable at the origin and has a pole in its
first derivative. Thus, the continuity module~\(\omega_{f_4}(\delta)\)
will decay very slowly for~\(\delta\to 0\), resulting in a slow
convergence of the interpolation process. 
\end{enumerate}

\subsection{Numerical Evaluation}
% If not stated otherwise, the interpolation is done for Leja ordered
% Chebyshev nodes.
Before the actual discussion of numerical results, let us briefly comment how
we evaluate discrete versions of the continuos-scale $L_1$, $L_2$
and $L_\infty$ norm, respectively.
To evaluate the numerical accuracy of the interpolation, instead of
computing the error norms directly, we approximate them numerically.
For the mentioned norms, we obtain as discrete analogues 
the~\(\mathfrak L_1\)- and~\(\mathfrak L_2\)-norms by the trapezoidal
rule, and the~\(\mathfrak L_\infty\)-norm by its discrete analogue on the
evaluation points for the approximation of the other norms. We employ
normalized norms on the test interval~\([-2,2]\), i.\,e.~\(\norm{1} =
1\) in any case. In order to capture the Gibbs phenomenon, we
integrate the~\(\mathfrak L_1\)- and~\(\mathfrak L_2\)-norms over~\(2n\)
strips.

%\subsection{Results}
%\label{sec:results}

  \begin{figure}
    \centering
    \includegraphics[width=\linewidth]{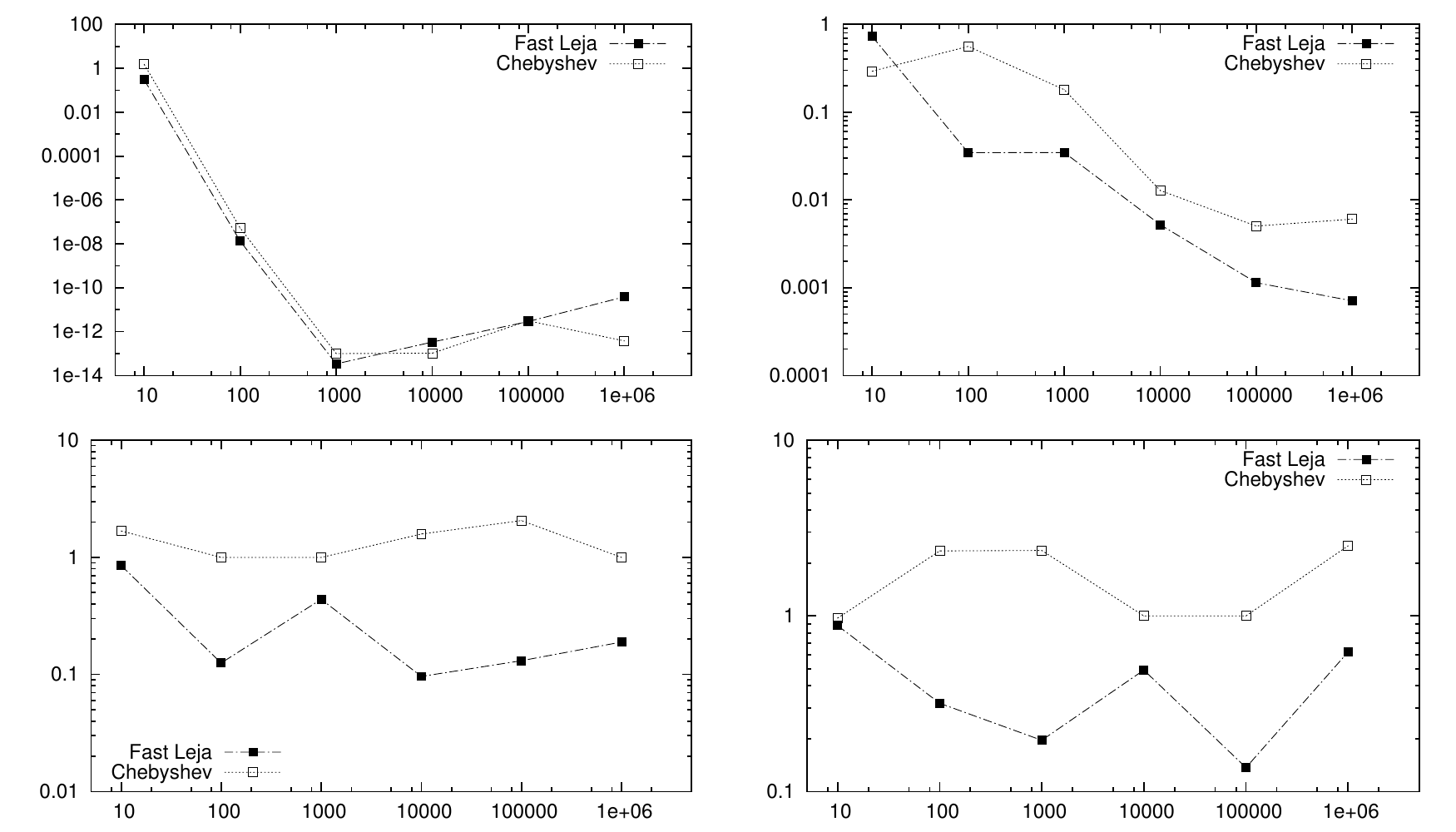}
    \caption{Comparison of approximation error for interpolation on
      Chebyshev and Fast Leja nodes,~\(\mathfrak L_\infty\)-norm:
      Runge function (top left), \(\sqrt{\abs{x}}\) (top right),
      Heaviside (bottom left), sawtooth (bottom right).}
    \label{fig:err-sup}
  \end{figure}
  \begin{figure}
    \centering
    \includegraphics[width=\linewidth]{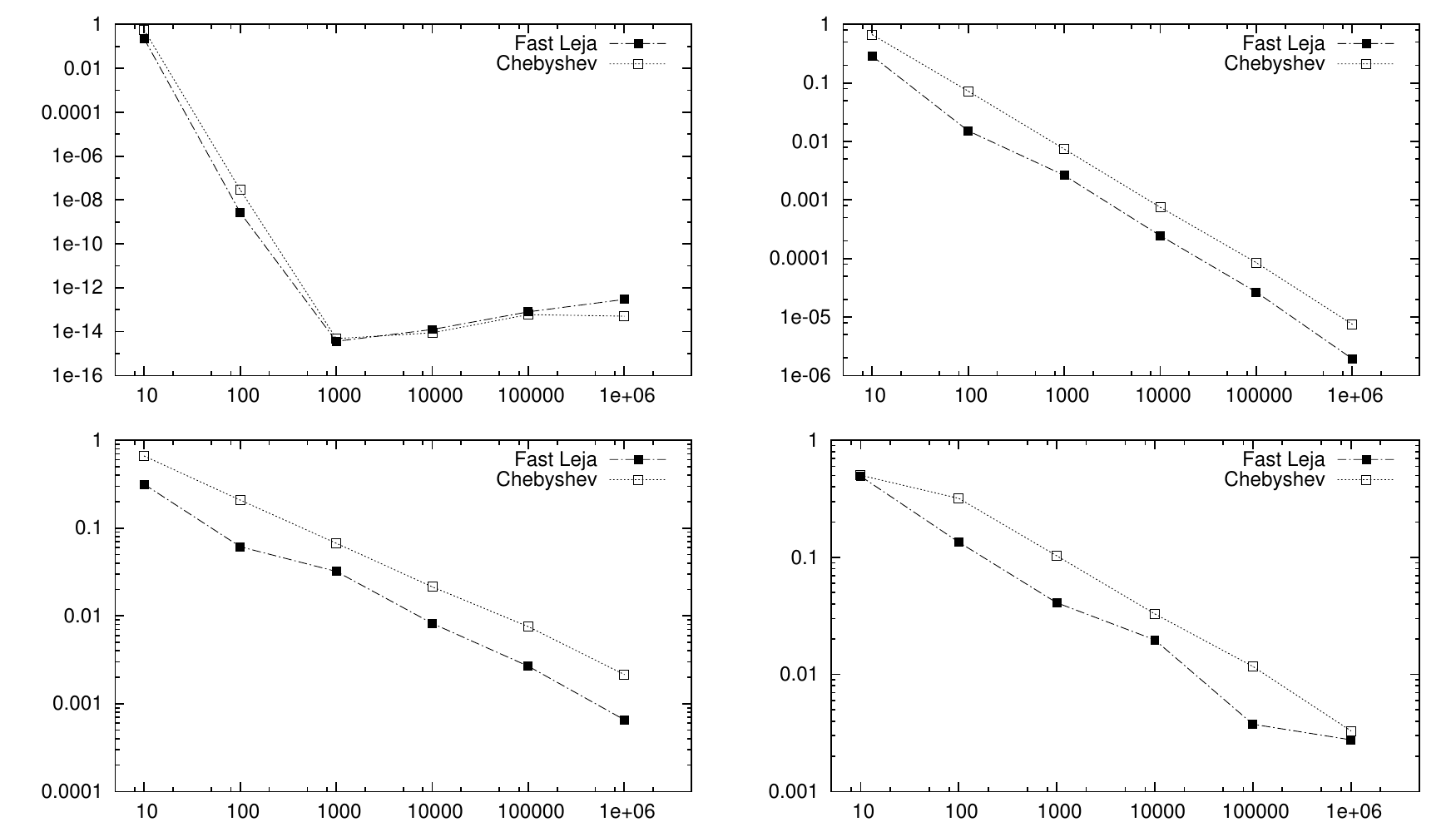}
    \caption{Comparison of approximation error for interpolation on
      Chebyshev and Fast Leja nodes,~\(\mathfrak L_2\)-norm:
      Runge function (top left), \(\sqrt{\abs{x}}\) (top right),
      Heaviside (bottom left), sawtooth (bottom right).}
    \label{fig:err-L2}
  \end{figure}
  \begin{figure}
    \centering
    \includegraphics[width=\linewidth]{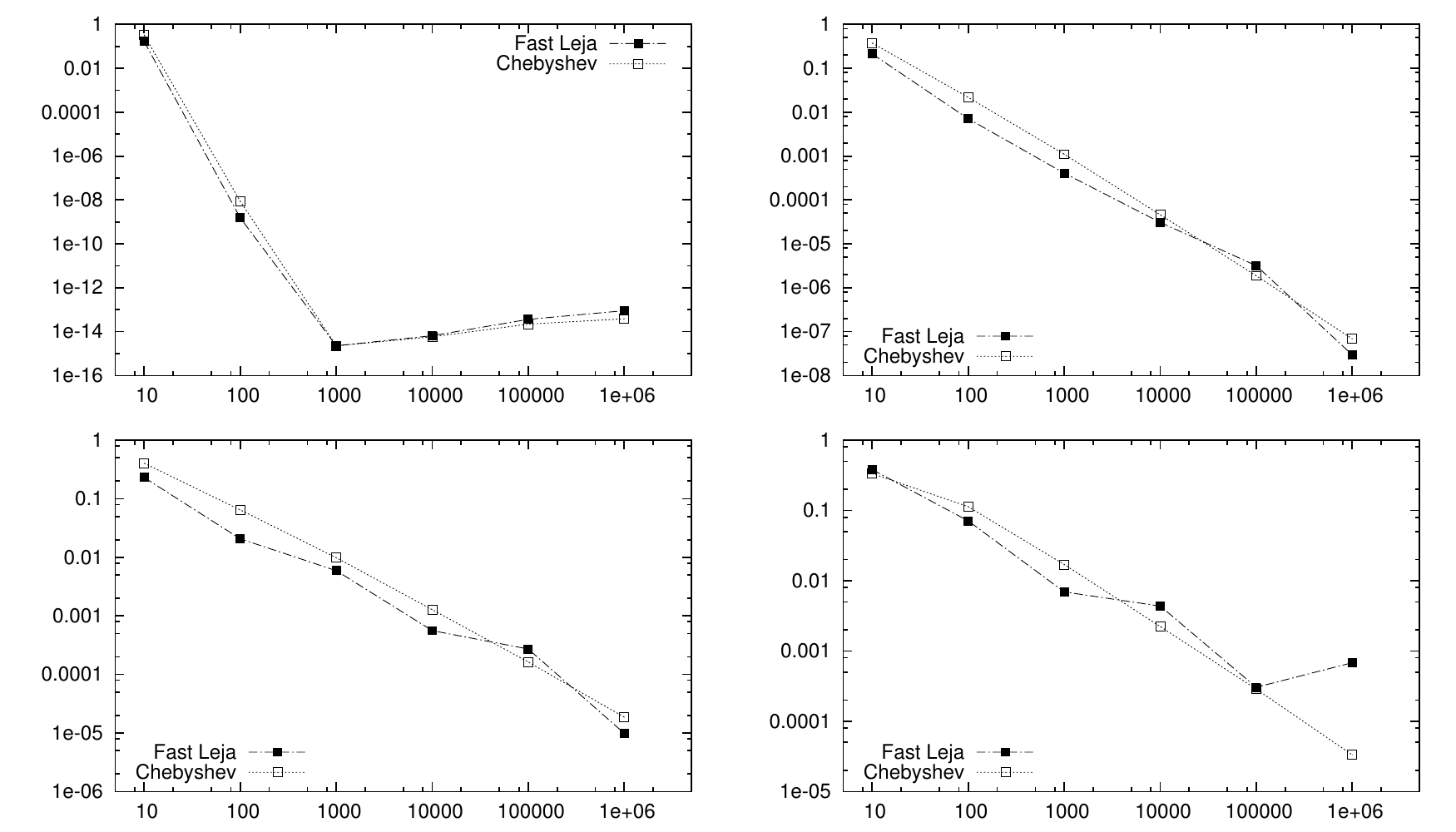}
    \caption{Comparison of approximation error for interpolation on
      Chebyshev and Fast Leja nodes,~\(\mathfrak L_1\)-norm:
      Runge function (top left), \(\sqrt{\abs{x}}\) (top right),
      Heaviside (bottom left), sawtooth (bottom right).}
    \label{fig:err-L1}
  \end{figure}

  \paragraph{General Observations on Numerical Convergence}
  In Figures~\ref{fig:err-sup}--\ref{fig:err-L1}, the approximation
  error for both the Newton interpolation on Leja ordered Chebyshev
  nodes and on Fast Leja nodes is presented. As predicted by the
  theory, for the interpolation of discontinuous functions there is no uniform
  convergence due to the Gibbs phenomenon. 
  This becomes particularly apparent when studying the \(\mathfrak L_\infty\)-norm
  (Figure~\ref{fig:err-sup}) of computational results. As a side note, we observe that FL points yield 
  slightly more accurate results in these cases. 
  However, the~\(\mathfrak L_1\)-norm (Figure~\ref{fig:err-L1}) 
  and even the~\(\mathfrak L_2\)-norm (Figure~\ref{fig:err-L2}) 
  become rather small. This corresponds to the fact that oscillations due to the Gibbs
  phenomenon do not decrease in size but become more localized at jump locations when
  increasing the degree of the interpolation polynomial. 
  For the classic example of the Runge function, which is~\(C^\infty\),
  even the~\(\mathfrak L_\infty\)-norm (Figure~\ref{fig:err-sup})
  of the error can be brought down to machine accuracy. 
 
  In the latter context, let us note that up to now in the literature, we did not recognize
  a study using comparatively high degrees of interpolation polynomials. For us it was at first
  amazing how stable the FL points perform having in mind their easy construction.

  \paragraph{Numerical Limit Distributions: Comparing Fast Leja Points and Chebyshev Points}
  The most important result of our investigation is that the error from the interpolation on Fast
  Leja points for high polynomial degrees approaches the error obtained on Leja ordered Chebyshev nodes, 
  note again that in the non-smooth cases it is even lower. Thus, there is no loss in accuracy, but still 
  much lower computational cost. Remember that the Fast Leja points can be
  computed beforehand for a standard interval and stored in a library. 
  The computational effort for the transformation to another interval is negligible.

  Let us elaborate here in addition on an important implication of our experiments. As already indicated, 
  it is a classic result that use of Chebyshev nodes is a pragmatic, computationally reasonable way 
  for high-order polynomial interpolation yielding a comparable quality as the best possible polynomial interpolation 
  \cite{natanson}, provided it can be resolved in a numerically stable way as done here via Leja ordering.
  As we have demonstrated one can achieve machine accuracy by going to extremely high numbers of interpolation 
  points. This means we have determined a regime of numerical convergence of the interpolation
  process. Therefore we have established in this regime also a reasonable numerical account of the 
  limit distribution of interpolation points. 

  Now let us consider the point that in this regime -- where we have a numerical account of the limit
  distribution of interpolation points -- the FL interpolation approaches the interpolation using Chebyshev points. 
  Note in this context that not only in this limit but also in the complete process of increasing the number 
  of interpolation points we always observe a comparable approximation quality, see the 
  Figures \ref{fig:err-sup} to \ref{fig:err-L1}.
  Because of the known favourable interpolation properties of the Chebyshev points which distinguish these
  from other choices, we have therefore validated numerically that the limit distribution of FL nodes tends to 
  the distribution of Chebyshev nodes.

  Let us add another comment to the last observation. It is quite clear now that we obtain a practical numerical analogon 
  to the theoretical assertion that the limit distribution of FL nodes tends to the distribution of Chebyshev nodes as 
  stated in \cite{baglama98}. However, beyond this assertion about the limit distribution which means to consider just the
  convergence region of very large number of points, we also observe in the whole process of increasing the number
  of interpolation points towards this convergence region a quantitatively similar behaviour of the interpolation
  with Fast Leja and Chebyshev points. This uniform behaviour of numerical accuracy that holds in addition to the behaviour in the
  numerical convergence zone is from our point of view an important observation that is up to now not reported in the 
  literature.

%%%%%%%%%%%%%%%%%%%%%%%%%%%%%%%%%%%%%%%%%%%%%%%%%5

\section{Conclusion}

Both the Leja ordering and the Fast Leja points enable impressive,
extremely high-degree Newton interpolation results with accuracy close
to machine accuracy. Especially for the Fast Leja points the
computational effort is rather small since the nodes might be read
from a library and just rescaled to the interpolation interval under
consideration. 

Fast Leja nodes enable the same interpolation quality as obtained
with Chebyshev nodes. At the same time they are a simple and computationally 
much more flexible tool. Therefore one may interpret the results of our study 
in the way that Newton polynomial interpolation with Fast Leja points is in practice 
the most efficient and simple way to do polynomial interpolation.
We think that the techniques of the Leja ordering and the Fast Leja points 
should be a topic mentioned in standard numerical analysis books.

\bibliographystyle{plain}      % mathematics and physical sciences

\bibliography{interp}   % name your BibTeX data base

\end{document}